\documentclass[12pt]{amsart}
\usepackage{amsmath, amsthm, amscd, amsfonts, amssymb, graphicx, color}
\newtheorem{thm}{Theorem}

\newtheorem{dfn}{Definition}

\begin{document}

\setcounter{page}{1}

\title[]{Furstenberg Theorem for Frequently Hypercyclic Operators}

\thanks{This work was supported by Seokyeong University in 2012. }
\maketitle

\begin{center}
\vspace{12mm} Eunsang Kim${}^{a,}$\footnote{eskim@hanyang.ac.kr}, \
and \; Tae Ryong Park
${}^{b,}$\footnote{trpark@skuniv.ac.kr} \

\vspace{5mm} ${}^a${\it
Department of Applied Mathematics ,\\ Hanyang University, Ansan
Kyunggi-do, Korea} \\
${}^b${\it
Department of Computer Engineering ,\\ Seokyeong University, Seoul, Korea}\\

\vspace{12mm}

\end{center}

\begin{center}
 {\bf Abstract}

\vspace{5mm}

\parbox{125mm}{In this paper, we  show that if the direct sum  $T\oplus T$ of frequently hypercyclic operators is frequently hypercyclic, then every higher direct sum $T\oplus\cdots\oplus T$ is also frequently hypercyclic.

} \vfill
\end{center}



\section{Introduction}

In this paper, we study the dynamics of linear operators on a separable $F$-space $X$. A bounded linear operator $T$ on $X$ is said to be hypercyclic if there is a vector $x\in X$ such that the orbit $O(x,T)=\{T^nx\mid n\in \mathbb{N}\}$ is dense in $X$. The operator $T$ is said to be topologically transitive if, for every pair of non-empty open subsets $U$ and $V$, there is an integer $n$ such that $T^nU\cap V\ne \emptyset$.  By the Baire category theorem, topological transitivity of $T$ is equivalent to the hypercyclicity of $T$. See \cite{book-GEandPeris} and  \cite{book-BayartandMatheron} for details and references. If $T\oplus T$ is hypercyclic, then the operator $T$ is called weakly mixing. It is shown in \cite{BesandPeris99} that the  weakly mixing property is equivalent to the Hypercyclicity Criterion. On the other hand, as shown in  \cite{RosaandRead} and  \cite{BayartandMatheron2},
hypercyclic operators may not be  weakly mixing, see also \cite{MatheronandBayart}. An interesting fact is so-called the Furstenberg theorem, which is given as follows: if $T$ is weakly mixing, then the $n$-fold product is $T\times\cdots\times T$ is weakly mixing for $n\ge 2$. The proof is given in \cite{book-GEandPeris} by using the 4-set trick. In the linear setting we have

\begin{thm}
Let $X$ be a separable $F$-space. If $T\oplus T$ is hypercyclic, then the higher sum $T\oplus\cdots \oplus T$ is also hypercyclic. \hfill$\qed$
\end{thm}

The $T$-orbit of a hypercyclic operator visits each non-empty open subsets of $X$. Then it is natural to ask  how often the orbit visits each non-empty open sets in $X$ and it leads to the notion of the frequently hypercyclic operators which has been introduced by
Bayart and Grivaux, see \cite{BayartandGrivaux} and \cite{BonillaandGrossErdmann}. In \cite{GrosseandPeris05}, it is shown that every frequently hypercyclic operator is weakly mixing. Based on ideas given in \cite{GrosseandPeris05}, we  prove the Furstenberg theorem for the frequently hypercyclic operators.

\section{Frequently Hypercyclic Operators}

Let $X$ be a separable $F$-space and  let $\mathcal{L}(X)$ be the space of continuous linear operators on $X$. By definition, an operator
$T\in \mathcal{L}(X)$ is hypercyclic if there is a vector $x\in X$ such that the orbit \[O(x,T)=\{T^nx\mid n\in\mathbb{N}\}\] is dense in $X$. In other words,  the  $T$-orbit $O(x,T)$ intersects with each non-empty open set $U$ in $X$. For a non-empty open subset $U$ of $X$, define
\[\mathbf{N}(x,U)=\{n\in\mathbb{N}\mid T^nx\in U\}.\]
If an operator $T$ on $X$ is hypercyclic, then there is a vector $x\in X$ such that for each non-empty open set $U$ in $X$, the set $\mathbf{N}(x,U)$ are all non-empty. For any non-empty open sets $U$ and $V$, let us define the {\it return set} as follows:
\[\mathbf{N}(U,V)=\{n\in \mathbb{N}\mid T^nU\cap V\ne\emptyset\}.\]
By the topological transitivit, if $T$ is hypercyclic, then each set $\mathbf{N}(U,V)$ is non-empty. 

If $T$ is weakly mixing, then there is a natural number $n$ such that for each open subsets $U_1, U_2$ and $V_1, V_2$ of $X$ such that
\[T^nU_1\cap V_1\ne\emptyset \text{ \ \ and \ \ } T^nU_2\cap V_2\ne\emptyset. \]
Then the $T\in \mathcal{L}(X)$ is weakly mixing if and only if 
\begin{align}\label{weakmixing}
\mathbf{N}(U_1,V_1)\cap \mathbf{N}(U_2,V_2)\ne\emptyset.
\end{align}
See \cite{GrosseandPeris05}, \cite{book-GEandPeris} and \cite{book-BayartandMatheron} for other formulas which are equivalent to (\ref{weakmixing}).


The frequently hypercyclicity corresponds to the largeness of each sets $\mathbf{N}(x,U)$, in other words, how frequently the $T$-orbit intersects with each open set $U$. Let us first recall that the lower density of  a subset $A$ in $\mathbb{N}$ which is given by
\[\underline{\text{dens}}(A)=\liminf_{N\to\infty}\frac{|A\cap [1,N]|}{N}\]
where $|A\cap [1,N]|$ denotes the cardinality of the set $A\cap [1,N]$.

\begin{dfn}{\rm Let $X$ be a topological vector space and let $T\in\mathcal{L}(X)$. The operator $T$ is called {\it frequently hypercyclic} if there is a vector $x\in X$ such that for every non-empty open set $U$, $\mathbf{N}(x,U)$ has positive lower density. Such a vector $x$ is called frequently hypercyclic for $T$ and the set of all frequently hypercyclic vectors for $T$ is denoted by $FHC(T)$.
}
\end{dfn}


If we enumerate an infinite set $A\subset\mathbb{N}$ as an increasing sequence $(n_k)_{k\in \mathbb{N}}$, then it is easy to see that $A$ has positive lower density if and only if there is a constant $C$ such that
\[n_k\le Ck \text{ \ \ \ \ for all  } k\ge 1\]
Thus, a vector $x\in X$ is frequently hypercyclic for $T$ if and only if for each non-empty open subset $U$ of $X$,   there is a strictly increasing sequence $(n_k)$  and some constant $C$ such that
\[T^{n_k}x\in U \text{ \ \ \ and \ \ \ } n_k\le Ck\]
for all $k\in \mathbb{N}$. We now prove the Furstenberg theorem for frequently hypercyclic operators.
\begin{thm}
Let $X$ be a separable $F$-space and let $T\in \mathcal{L}(X)$. If $T\oplus T$ is frequently hypercyclic, then $3$-fold sum $T\oplus T \oplus T$ is also frequently hypercyclic.
\end{thm}

{\it Proof. } We will show that there is a vector $x_1\oplus x_2\oplus x_3\in X\oplus X\oplus X$ satisfying the following property: for each non-empty open subsets $U_1$, $U_2$ and $U_3$ of $X$, there is a  strictly increasing sequence  $(n_k)_{k\in \mathbb{N}}$ of natural numbers and a constant $C$ such that for $i=1,2,3$,
\[T^{n_k}x_i\in U_i \text{ \ \  and  \ \ } n_k\le Ck \text{ \ \ \   for all }k\in  \mathbb{N}\]

First, we note that if $T\oplus T$ is  frequently hypercyclic, then $T\oplus T$ is  hypercyclic. By the Furstenberg theorem $T\oplus T \oplus T$ is also hypercyclic. Thus there is a hypercyclic vector  $x_1\oplus x_2\oplus x_3\in X\oplus X\oplus X$ such that for each non-empty open subsets $U_1$, $U_2$ and $U_3$ of $X$,
\[\mathbf{N}(x_1,U_1)\cap \mathbf{N}(x_2,U_2)\cap \mathbf{N}(x_3,U_3)\ne\emptyset\]

Suppose that $x_1\oplus x_2\in FHC(T\oplus T)$. Then for non-empty open sets $U_1$ and $U_2$ there is a strictly increasing sequence $(m_k)_{k\in\mathbb{N}}$ and a constant $C_1$ such that
\[T^{m_k}x_1\in U_1, \ \ T^{m_k}x_2\in U_2 \text{ \ \ \ and \ \ } m_k\le C_1k\]
Since $T$ is hypercyclic, the set $\mathbf{N}(U_1,U_2):=\{l\in \mathbb{N}\mid T^lU_1\cap U_2\ne\emptyset\}$ is non-empty and the $T$-orbit $O(x_1,T)$ is dense in $X$, there is an increasing sequence $(b_j)_{j\in \mathbb{N}}$ such that
\begin{align}\label{bn}
x_2=\lim_{j\to\infty}T^{b_j}x_1
\end{align}
Since $T$ is continuous,
\begin{align}
T^{m_k}x_2=\lim_{j\to\infty}T^{b_j}T^{m_k}x_1\in U_2
\end{align}
Thus there is an integer $N$ such that for all $j\ge 1$
\[T^{b_{j+N}}T^{m_k}x_1\in U_2\]
In particular, for all $k\in \mathbb{N}$
\[T^{b_{k+N}}T^{m_k}\in U_2 \text{ \ \ and  \ \ } T^{b_{k+N}}U_1\cap U_2\ne\emptyset \]
In other words, the sequence $(b_{k+N})_{k\in\mathbb{N}}$ is in $\mathbf{N}(U_1,U_2)$ and since $x_1$ and $x_2$ are frequently hypercyclic, the sequence $(b_{k+N})_{k\in\mathbb{N}}$ satisfies $b_{k+N}=O(k)$ ({\it cf.} \cite{MatheronandBayart}). Let
\[U_{1k}=U_1\cap T^{-b_{k+N}}U_2\]
Then $\mathbf{N}(x_1,U_{1k})\subset \mathbf{N}(x_1,U_1)$.  If $l\in\mathbf{N}(x_1,U_{1k})$ then for all $k\in \mathbb{N}$,
\[T^lx_1\in U_1 \text{ \ \ and \ \ } T^{b_{k+N}}T^l x_1\in U_2\]
By (\ref{bn}) and $T^{b_{k+N}}T^l x_1=T^lT^{b_{k+N}}x_1$, we get
$l\in \mathbf{N}(x_1, U_{1k})$. Thus
\begin{align}\label{con1}
\mathbf{N}(x_1,U_{1k})\subset \mathbf{N}(x_1,U_{1})\cap \mathbf{N}(x_2,U_{2})
\end{align}

Applying the same argument for $\mathbf{N}(x_2,U_2)\cap \mathbf{N}(x_3,U_3)$ we
may obtain
\begin{align}\label{con2}
\mathbf{N}(x_2,U_{2j})\subset \mathbf{N}(x_2,U_{2})\cap \mathbf{N}(x_3,U_{3})
\end{align}
Now by (\ref{con1}) and (\ref{con2})
\begin{align}
\mathbf{N}(x_1,U_{1k})\cap \mathbf{N}(x_2,U_{2j})\subset \mathbf{N}(x_1,U_1)\cap \mathbf{N}(x_2,U_2)\cap \mathbf{N}(x_3,U_3)
\end{align}

Since $x_1\oplus x_2\in FHC(T\oplus T)$, there is a strictly increasing sequence $(n_k)_{k\in \mathbb{N}}$, which may be an enumeration of thet set $\mathbf{N}(x_1,U_{1k})\cap \mathbf{N}(x_2,U_{2j})$, such that for some constant $C$, and  for $i=1,2,3$,
\[T^{n_k}x_i\in U_i \text{ \ \  and  \ \ } n_k\le Ck \text{ \ \ \   for all }k\in  \mathbb{N}\]
as desired. \hfill\qed

By proceeding induction, we have the main result
\begin{thm}
Let $X$ be a separable $F$-space and let $T\in \mathcal{L}(X)$. If $T\oplus T$ is frequently hypercyclic, then the higher product $T\oplus \cdots\oplus T$ is also frequently hypercyclic.\hfill\qed
\end{thm}


\end{document}